# Cubic Crystals and Vertex-Colorings of the Cubic Honeycomb

Mark L. Loyola[1], Ma. Louise Antonette N. De Las Peñas[1], Antonio M. Basilio[2]
[1]Mathematics Department, [2]Chemistry Department
Ateneo de Manila University
Loyola Heights, Quezon City 1108, Philippines
mloyola@math.admu.edu.ph, mlp@math.admu.edu.ph, abasilio@ateneo.edu

**Abstract**

In this work, we discuss vertex-colorings of the cubic honeycomb and we illustrate how these colorings can demonstrate the structure and symmetries of certain cubic crystals.

## Introduction

In modern crystallography, a real crystal is modeled as a 3-dimensional array of points (a lattice) with each point decorated by a motif, which usually stands for an atom, a molecule or an aggregate of atoms, that is being repeated or translated periodically throughout the crystal model (**Figure 1**). Alternative ways of modeling crystals include the use of space-filling polyhedra, sphere-packings, infinite graphs, and regular system of points [5].

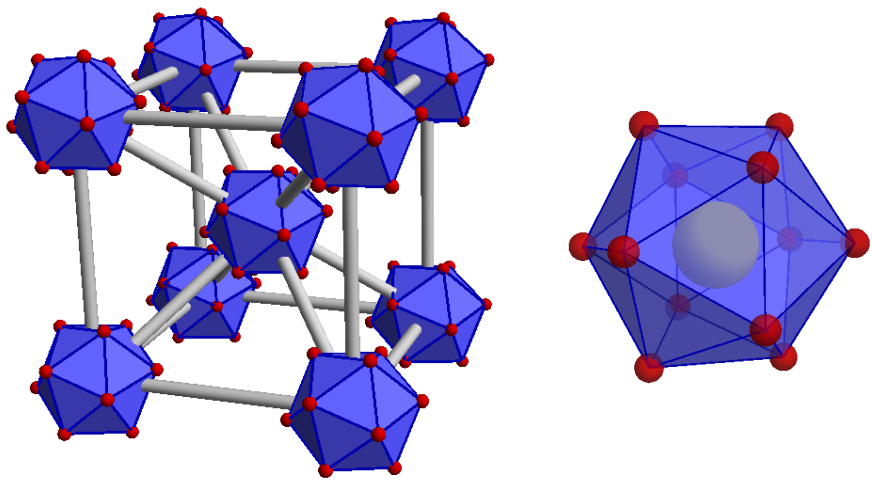

**Figure 1:** *The $Al_{12}W$ crystal modeled as a body-centered cubic lattice decorated at each point with an icosahedron consisting of an aluminum (Al) atom at each vertex and a tungsten (W) atom at the center.*

The models mentioned above exemplify two important physical properties or aspects of crystals – their high symmetry and their periodic nature. Yet another approach of representing crystals that demonstrates these two properties is through the use of colorings of periodic structures such as lattices or tessellations. Colors in this context of coloring can signify any meaningful attribute of a crystal such as the different types of its discrete units (atoms, molecules, group of atoms) or the magnetic spins associated with its atoms. Colors in colored models of crystalline structures are particularly useful in describing, visualizing, and analyzing the ordered arrangement and relative positions of the discrete units that gives a crystal its shape and structure.

This paper describes a systematic way of representing certain cubic crystals via vertex-colorings of the cubic honeycomb (tessellation of the space by cubes) where a color in the coloring signifies a particular type of atom. These colorings are obtained by computing the orbits of the vertices of the cubic honeycomb under the group action of a subgroup $H$ of the symmetry group $G$ of the uncolored honeycomb where we require that the elements of $H$ induce a permutation of the colors. Each $H$-orbit is then partitioned using the left cosets of a subgroup $J$ of $H$.

## Symmetry Group of the Cubic Honeycomb

The *cubic honeycomb* $\{4, 3, 4\}$ shown in **Figure 3(a)** is defined as the regular space-filling tessellation of the Euclidean space $\mathbf{E}^3$ made up of cubes or *cubic cells*. The cubes are arranged face to face in space so that each edge in the tessellation is shared by exactly four cubic cells. It belongs to the family of 28 convex uniform honeycombs [3] and is also characterized as the 3-dimensional analogue of the $\{4, 4\}$ regular tiling of the plane by squares.

The isometries of $\mathbf{E}^3$ that send the cubic honeycomb $\{4, 3, 4\}$ to itself form what we refer to as the symmetry group $G$ of the honeycomb. A set of generators for $G$ consists of the reflections $P, Q, R, S$ about the planes $m_P, m_Q, m_R, m_S$. A cubic cell presented in **Figure 2** from the $\{4, 3, 4\}$ honeycomb in **Figure 3(a)** shows these planes of reflections. The plane $m_P$ is parallel to the top and bottom faces of the cube and cuts the other faces in half; the plane $m_Q$ is the diagonal plane which cuts the front and back faces of the cube into right triangular regions; similarly, the plane $m_R$ is also a diagonal plane, this time, cutting the top and bottom faces into right triangular regions; finally, the plane $m_S$ is the plane which is off the center of the cube and contains the front face. These four planes form the walls of a tetrahedron with dihedral angles $\pi/4, \pi/3, \pi/4, \pi/2, \pi/2, \pi/2$, which serves as a fundamental polyhedron of the tessellation.

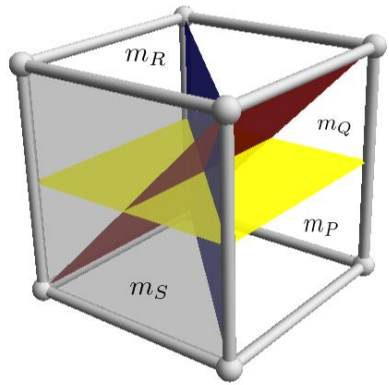

**Figure 2**: *The planes $m_P, m_Q, m_R, m_S$ of the reflections $P, Q, R, S$ of the symmetry group of the cubic honeycomb viewed on a cubic cell.*

As a reflection group $G$ is isomorphic to the affine Weyl group $C_3$ with Coxeter presentation

$$< P, Q, R, S \mid P^2 = Q^2 = R^2 = S^2 = (PQ)^4 = (QR)^3 = (RS)^4 = (PR)^2 = (PS)^2 = (QS)^2 = 1 >$$

and corresponding Coxeter diagram 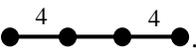.

The cubic honeycomb is a *vertex-transitive* (*isogonal*) tessellation, that is, for each pair of vertices $x_1$ and $x_2$ in the honeycomb, we can find a symmetry $g \in G$ for which $gx_1 = x_2$. Thus, each vertex of the honeycomb may be sent to any other vertex by some symmetry in $G$. The notion of vertex transitivity is important if we are to consider vertex-colorings of the cubic honeycomb which will be apparent from the discussion in the next section.

## A Method to Obtain Vertex-Colorings of the Cubic Honeycomb

A general framework for obtaining symmetrical patterns [2] will be adopted to arrive at vertex-colorings of the {4, 3, 4} honeycomb, which we state as follows:

Let $G$ be the symmetry group of the uncolored {4, 3, 4} honeycomb and $X$ the set of vertices in the honeycomb. If $C = \{c_1, c_2, \ldots, c_n\}$ is a set of $n$ colors, an onto function $f : X \to C$ is called an *n-coloring* of $X$. Each $x \in X$ is assigned a color in $C$. The coloring determines a partition $\Pi = \{f^{-1}(c_i) : c_i \in C\}$ where $f^{-1}(c_i)$ is the set of elements of $X$ assigned color $c_i$.

Let $H$ be the subgroup of $G$ which consists of symmetries in $G$ that effect a permutation of the colors in $C$. Then $h \in H$ if for every $c \in C$, there is a $d \in C$ such that $h(f^{-1}(c)) = f^{-1}(d)$. This defines an action of $H$ on $C$ where we write $hc := d$ if and only if $h(f^{-1}(c)) = f^{-1}(d)$. Since $H$ acts on $C$, there exists a homomorphism $\sigma$ from $H$ to Perm($C$) where Perm($C$) is the group of permutations of $C$.

Now denote by $Hc_i$ the orbit of the color $c_i$ under the action of $H$ on $C$. Let $J_i = \{h \in H : hc_i = c_i\}$ be the stabilizer of $c_i$ in $H$. Pick an element from each $H$-orbit of $X$ with an element colored $c_i$. Put these elements together in a set $X_i$. Then the set of all elements of $X$ that are colored $c_i$ is $J_iX_i = \{jx : j \in J_i, x \in X_i\}$, that is, $f^{-1}(c_i) = J_iX_i$. A one-to-one correspondence results between the sets $Hc_i$ and $\{hJ_iX_i : h \in H\}$ where $hJ_iX_i$ denotes the image of $J_iX_i$ under $h$. These assumptions lead to the following result:

**Theorem:**

1. The action of $H$ on $Hc_i$ is equivalent to its action on $\{hJ_i : h \in H\}$ by left multiplication.
2. The number of colors in $Hc_i$ is equal to $[H : J_i]$.
3. The number of $H$-orbits of colors is at most the number of $H$-orbits of elements of $X$.
4. If $x \in X_i$ and Stab$_H(x) = \{h \in H : hx = x\}$ is the stabilizer of $x$ under the action of $H$ on $X$ then
   (a) Stab$_H(x) \leq J_i$ .
   (b) $|Hx| = [H : J_i] \cdot [J_i : $ Stab$_H(x)]$

We now give the process to obtain vertex-colorings of the {4, 3, 4} honeycomb where the elements of $H$ effect a permutation of the colors in the given coloring.

1. Determine the different $H$-orbits of vertices in $X$ and color each $H$-orbit separately.
2. To color an $H$-orbit $Hx$, $x \in X$, choose a subgroup $J$ of $H$ such that Stab$_H(x) \leq J$.
3. Assign a color to every element of the set $\{hJx: h \in H\}$. If $[H : J] = k$, then $Jx$ is $1/k$ of the vertices in the $H$-orbit where $x$ belongs. The set $Jx$ is given one color and each of the remaining $k - 1$ elements of the set gets a different color.
4. If two $H$-orbits of vertices are to have a color in common, the subgroup $J$ used should contain the stabilizers of representative vertices from the two $H$-orbits.

To illustrate the method of obtaining vertex-colorings of the {4, 3, 4} honeycomb, we consider the following examples:

**Illustration 1.** As a first example, we construct a vertex-coloring of the {4,3,4} honeycomb where all the elements of its symmetry group $G = <P, Q, R, S>$ effect a permutation of the colors. We start with the vertex labeled $x$ in **Figure 3(b)**. Note that the vertices of the honeycomb form one $G$-orbit of vertices, that is $Gx = X$. The stabilizer of $x$ in $G$, $Stab_G(x)$, is the group of 48 elements generated by the reflections $Q$, $R$, $PQRSRQP$, a group of type $O_h$, also known as the *octahedral group*. To color $X$, we select a subgroup $J$ of $G$ satisfying the condition that $Stab_G(x) \leq J$. With the aid of GAP [6], we choose the index 2 subgroup $J = <Q, R, S, PQP>$ of $G$. We assign to $Jx$ the color light blue. To color the rest of $X$, we apply the reflection $P$ about the plane $m_P$ to $Jx$ to obtain the 2-coloring shown in **Figure 3(b)**.

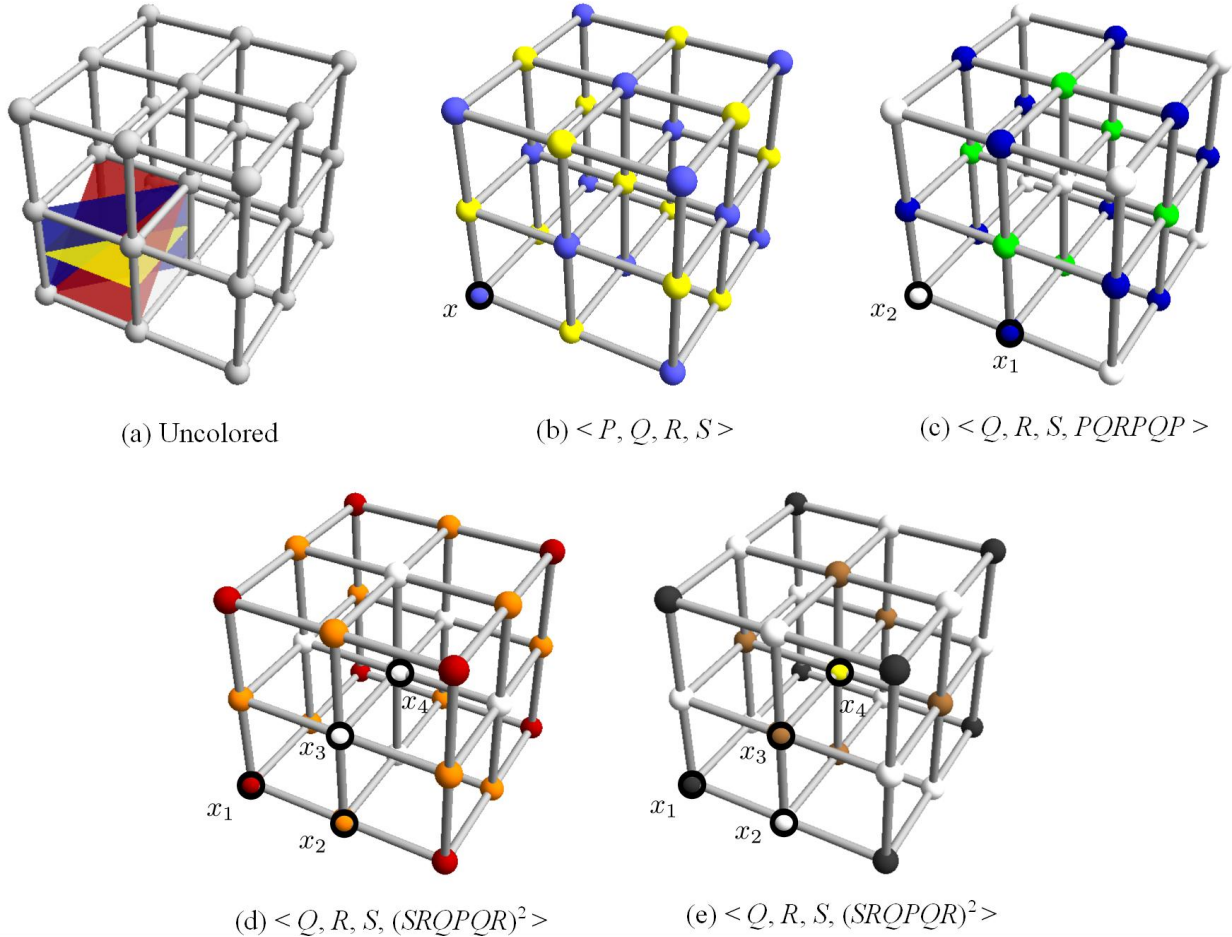

**Figure 3**: (a) uncolored honeycomb with planes of reflections (generators) for its symmetry group $G$; (b) – (e): *Vertex-colorings of the cubic honeycomb* {4, 3, 4}. *The subgroup $H$ of $G$ consisting of elements that effect a permutation of the colors in each coloring is indicated. A representative vertex of each $H$-orbit is labeled $x_i$.*

**Illustration 2.** In this example, we construct a vertex-coloring of the {4, 3, 4} honeycomb where the elements of the group $H = <Q, R, S, PQRPQP>$ effect a permutation of the colors in the given coloring. In this case, $H$ is an index 4 subgroup of $G = <P, Q, R, S>$. Upon inspection, we find that there are two $H$-orbits of $X$, we call these orbits $X_1 = Hx_1$ and $X_2 = Hx_2$, where $x_1$, $x_2$ are representative atoms coming from the two $H$-orbits. One way to color $X_1$ is to use an index 2 subgroup $J_1 = <Q, R, S, (SRQPQR)^2>$ of $H$ that contains the $Stab_H(x_1) = <Q, S, (QPQRQPQS)^2>$. To color $X_2$, we choose $H$, which contains $Stab_H(x_2) = <Q, R, PQRSRQP>$. Assigning the color dark blue to $J_1x_1$, green to $(PQRPQP)Jx_1$ in orbit $X_1$ and white to $X_2 = Hx_2$ will result in the 3-coloring given in **Figure 3(c)**.

Other non-perfect colorings of the {4, 3, 4} honeycomb are presented in **Figures 3(d)** and **3(e).** For these colorings the elements of $H_2 = <Q, R, S, (SRQPQR)^2>$ consists of elements that permute the colors. There are four $H_2$-orbits formed namely, $H_2 x_1$, $H_2 x_2$, $H_2 x_3$ and $H_2 x_4$. The coloring in Figure 3(d) is obtained by assigning the color red to $H_2 x_1$, orange to $H_2 x_2$ and white to $H_2 x_3 \cup H_2 x_4$ (two $H$-orbits share a color). For the coloring shown in Figure 3(e), black is assigned to $H_2 x_1$, white to $H_2 x_2$, brown to $H_2 x_3$ and yellow to $H_2 x_4$.

## Cubic Crystals from Vertex-Colorings of the Cubic Honeycomb

We now discuss a manner of representing certain crystals with cubic structures [1] via vertex-colorings of the cubic honeycomb {4, 3, 4}. In the representation scheme, a type of atom in the corresponding cubic crystal is represented by a particular color in the vertex-coloring.

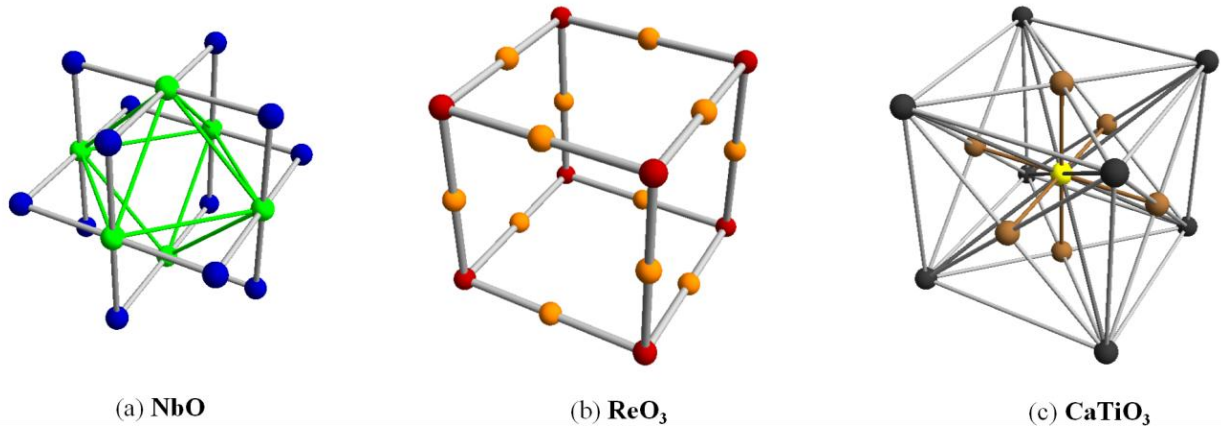

(a) **NbO**  (b) **ReO₃**  (c) **CaTiO₃**

**Figure 4**: *Unit cells of the cubic crystals* **NbO**, **ReO₃** *and* **CaTiO₃**.

The 2-coloring discussed in **Illustration 1** and presented in **Figure 3(b)** is an example of a *perfect* coloring of the {4, 3, 4} honeycomb. By a perfect coloring we mean all the elements of the symmetry group $G$ of the uncolored honeycomb effect a permutation of colors of the given coloring. In any perfectly colored pattern, the colors are distributed evenly throughout the coloring. Because of this property, a perfect 2-coloring of the {4, 3, 4} honeycomb may be used to model a cubic crystal with an even distribution of two atoms, such as the salt crystal **NaCl**. Among the cubic crystals that can be modeled by this 2-coloring include silver chloride (**AgCl**), calcium oxide (**CaO**), sodium fluoride (**NaF**), and potassium bromide (**KBr**).

The coloring given in **Illustration 2** and presented in **Figure 3(c)** is a *non-perfect* coloring of the {4, 3, 4} honeycomb. There is an uneven distribution of colors in the given coloring and not all elements of $G$ effect a permutation of the given colors. One can verify that all the symmetries of the group $H = <Q, R, S, PQRPQP>$ permute the given colors. The coloring can be used to model the crystal niobium monoxide (**NbO**) with unit cell shown in **Figure 4(a).** The niobium (**Nb**) atoms are represented by the green vertices while the oxygen (**O**) atoms are represented by the dark blue vertices. In this representation, the vertices colored white do not represent a particular type of atom. The **Nb** and **O** atoms belong to the same orbit under the reflections $Q$, $R$, the 2-fold rotation $PQRPQP$, as well as other elements of $H$.

The coloring appearing in **Figure 3(d)** can be used to represent a cubic crystal such as rhenium trioxide (**ReO₃**) where the rhenium (**Re**) atoms are represented by the red vertices and the oxygen atoms are represented by orange vertices. The distribution of atoms in the orbits $H_2 x_1$, $H_2 x_2$ occur in the ratio 1:3. As in the previous example, the vertices colored white do not represent any type of atom. Other cubic crystals that can be modeled using this coloring are copper nitride ($Cu_3N$) and tungsten trioxide ($O_3W$).

Finally, the coloring appearing in **Figure 3(e)** may correspond to a cubic crystal such as perovskite ($CaTiO_3$). The calcium (**Ca**), titanium (**Ti**), and oxygen atoms are represented by the black, yellow and brown vertices, respectively. The distribution of atoms in the orbits $H_2 x_1$, $H_2 x_4$ and $H_2 x_3$ occur in the ratio 1:1:3. Once again, the white vertices do not represent any atom. The other cubic crystals that can be represented by this coloring are barium titanate ($BaTiO_3$), lead zirconate ($PbZrO_3$), and lead titanate ($PbTiO_3$).

## Conclusion

In this work, we have highlighted the use of vertex-colorings of the cubic honeycomb to describe the symmetrical structure of cubic crystals. In general, in arriving at a coloring of the set *X* of vertices of the {4, 3, 4} cubic honeycomb that will model a given cubic crystal, an important point of consideration would be to determine the group *H* consisting of symmetries that effect a permutation of the atoms in the given crystal. The distribution of atoms in the given crystal will suggest the choice of subgroups to be used in arriving at colorings of the *H*-orbits of *X*.

Color symmetry and group theory are important tools to facilitate the understanding of crystal structures. By investigating and constructing vertex-colorings of honeycombs or lattices, some properties of crystals with various types of atoms may be predicted based on their symmetrical arrangements. This is not included in the scope of this paper and can be considered for future work.